\documentclass[12pt,a4paper,twoside,english,american]{scrartcl}
\usepackage{helvet}
\usepackage{babel}
\usepackage{verbatim}
\usepackage{uni_math}
\usepackage{amsthm}
\usepackage{amsmath}
\usepackage{amssymb}
\usepackage[unicode=true,pdfusetitle,
 bookmarks=true,bookmarksnumbered=false,bookmarksopen=false,
 breaklinks=false,pdfborder={0 0 1},backref=false,colorlinks=false]
 {hyperref}

\makeatletter

\numberwithin{equation}{section}
\numberwithin{figure}{section}
 \theoremstyle{definition}
 \newtheorem*{defn*}{\protect\definitionname}
\theoremstyle{plain}

  \theoremstyle{remark}
  \theoremstyle{plain}
  
  \theoremstyle{plain}
  
  \theoremstyle{plain}
  
  \theoremstyle{definition}

\usepackage[T1]{fontenc}    

\usepackage{a4wide}        
\usepackage{tu-preprint}
\usepackage{amsmath}
\usepackage{enumerate}

\usepackage{euscript}
\usepackage{mathtools}
\allowdisplaybreaks[1] 
\usepackage{amsfonts}
\newenvironment{keywords}{ \noindent\footnotesize\textbf{Keywords and phrases:}}{}

\newenvironment{class}{\noindent\footnotesize\textbf{Mathematics subject classification 2010:}}{}

\theoremstyle{definition}
\newtheorem{hyp}[subsection]{Hypotheses}

\usepackage{color,tu-preprint}

\DeclareMathAccent{\Circ}{\mathalpha}{operators}{"17}

\renewcommand{\Re}{\operatorname{\mathfrak{Re}}}

\renewcommand{\tilde}{\widetilde}
\renewcommand*{\epsilon}{\varepsilon}
\renewcommand*{\theta}{\vartheta}
\renewcommand*{\rho}{\varrho}

\makeatother

  \addto\captionsamerican{\renewcommand{\corollaryname}{Corollary}}
  \addto\captionsamerican{\renewcommand{\definitionname}{Definition}}
  \addto\captionsamerican{\renewcommand{\lemmaname}{Lemma}}
  \addto\captionsamerican{\renewcommand{\propositionname}{Proposition}}
  \addto\captionsamerican{}
  \addto\captionsamerican{\renewcommand{\theoremname}{Theorem}}
  \addto\captionsenglish{\renewcommand{\corollaryname}{Corollary}}
  \addto\captionsenglish{\renewcommand{\definitionname}{Definition}}
  \addto\captionsenglish{\renewcommand{\lemmaname}{Lemma}}
  \addto\captionsenglish{\renewcommand{\propositionname}{Proposition}}
  \addto\captionsenglish{}
  \addto\captionsenglish{\renewcommand{\theoremname}{Theorem}}
  \providecommand{\corollaryname}{Corollary}
  \providecommand{\definitionname}{Definition}
  \providecommand{\lemmaname}{Lemma}
  \providecommand{\propositionname}{Proposition}
  
\providecommand{\theoremname}{Theorem}

\begin{document}
\selectlanguage{english}%
\institut{Institut f\"ur Analysis}

\preprintnumber{MATH-AN-08-2013}

\preprinttitle{On Non-Autonomous Integro-Differential-Algebraic Evolutionary Problems.}

\author{Marcus Waurick} 

\makepreprinttitlepage

\selectlanguage{american}%
\setcounter{section}{-1}

\date{}

\title{On Non-Autonomous Integro-Differential-Algebraic Evolutionary Problems.}

\author{Marcus Waurick\\
Institut f\"ur Analysis, Fachrichtung Mathematik\\
Technische Universit\"at Dresden\\
Germany\\
marcus.waurick@tu-dresden.de\\
}
\maketitle
\begin{abstract} In this article, we show that a technique for showing well-posedness results for evolutionary equations in the sense of \cite{Picard} established in \cite{RainerPicard2013} applies to a broader class of non-autonomous integro-differential-algebraic equations. Using the concept of evolutionary mappings we prove that the respective solution operators do not depend on certain parameters describing the underlying spaces in which the well-posedness results are established.
\end{abstract}
\begin{keywords}
non-autonomous, evolutionary problems, extrapolation spaces, in\-te\-gro-differen\-tial-al\-ge\-braic equations \end{keywords}

\begin{class}
35F05, 35F46, 37L05, 35M10
\end{class}

\newpage

\tableofcontents{}


\cleardoublestandardpage

\section{Introduction}

It has been investigated that linear (evolutionary) equations of mathematical physics share a common form, see e.g.~\cite{Picard}. In a Hilbert space framework, formally, these equations can be written as
\begin{equation}\label{eq:0}
   \dot v (t,x) + A u(t,x) = f(t,x),
\end{equation}
for $t\in [0,\infty)$ and $x$ belonging to a certain open subset of $\R^n$. The notation $\dot v$ stands for the time-derivative of $v$. The operator $A$ is a maximal monotone ($m$-accretive) linear operator containing the spatial derivatives and $f$ is a certain source term. As the functions $u$ and $v$ are both unknowns, equation \eqref{eq:0} is not enough to uniquely determine $u$ and $v$. Thus, equation \eqref{eq:0} needs to be supplemented by the so-called ``material law'' or ``constitutive relation''. The material law links $u$ and $v$ via a bounded linear operator $\s M$ acting in space-time in the way that
\begin{equation}\label{eq:1}
  v= \s M u.
\end{equation}
Of course equations \eqref{eq:0} and \eqref{eq:1} are subject to certain initial conditions, which -- for simplicity --  we assume to be $0$. Boundary conditions are encoded in the domain of $A$. 

Plugging \eqref{eq:1} into \eqref{eq:0}, we arrive at
\begin{equation}\label{eq:2}
    \left(\s M u\right)^\cdot  + A u = f,
\end{equation}
In this article, we discuss well-posedness issues of the equation \eqref{eq:2} in a particular Hilbert space setting. We shall note here that, if $\s M = 1$, equation \eqref{eq:2} is well-known to be well-posed in the sense that $-A$ generates a $C_0$-semigroup and the respective solution can be obtained via the variation of constants formula, see e.g.~\cite{EngNag}. By suitably weighting the norm and assuming that $\s M$ only acts on the spatial variables and -- at the same time -- is continuously invertible, we also realize that equation \eqref{eq:2} can be discussed within a semi-group framework. Writing $\partial_0$ for time-differentiation ($\partial_0^{-1}$ is then integration, which is made more precise later on) and assume that $\s M$ is given in the form $M+\partial_0^{-1}N$ for suitable bounded linear operators $M$ and $N$, we infer that equation \eqref{eq:2} reads as
\begin{equation}\label{eq:3}
  \partial_0\left( \s M u\right)+A u = \partial_0\left( M + \partial_0^{-1} N\right)u + A u = \partial_0 M u + Nu + Au = f.
\end{equation}
Of course, if $M$ is continuously invertible, it might be possible to show well-posedness of \eqref{eq:3} with semi-group techniques by suitable perturbation theorems. However, if $M$ has a non-trivial nullspace the equation \eqref{eq:3} amounts to be a differential-algebraic system, which might be hard to treat within a $C_0$-semi-group framework. Moreover, for example if $M=0$, (time-)regularity of a possible solution of \eqref{eq:3} cannot be expected. 

Since we assume $\s M$ to be non-autonomous in general, we need to incorporate other solution techniques. The strategies discussed in \cite{Acquistapace1987,Kato1953,Tanabe1997} generalize the semi-group perspective in the sense that the generators are no longer time independent. In order to apply the theory developed in \cite{Acquistapace1987,Kato1953,Tanabe1997},  the operator $M$ again has to be inverted. The authors of \cite{Acquistapace1987} focus on the parabolic case.

 Thus, other techniques need to be incorporated. Realizing that \eqref{eq:3} is a formal sum of the two (unbounded) operators $\partial_0M+N$ and $A$, we seek for conditions, which guarantee that this operator sum is closable and continuously invertible.  The authors of \cite{Prato1975} give general conditions, which guarantee that a sum of two unbounded linear operators is closable and (that its closure is) continuously invertible. They mainly focus on the parabolic case, i.e., they assume that one of the operators involved is sectorial. However, in \cite[Section 5]{Prato1975} they also discuss the hyperbolic case in a Banach space setting. In order to derive closability  and surjectivity of the operator sum they assume certain resolvent estimates and properties of how the Yosida approximants converge. In \cite{RainerPicard2013}, such estimates are not assumed. However, in \cite{RainerPicard2013} one sticks to the Hilbert space case, which makes things conceptually easier. In \cite{Showalter1974} a related problem class has been studied. However, since the operators under consideration are defined via forms, the author focussed on the parabolic/elliptic case.

The heart of the solution theory in the present article, and thus the possibility to define the operator sum and to derive that the respective operator is continuously invertible, is a positive definiteness constraint in space-time. This allows for the consideration of equations with change of type switching from elliptic, to hyperbolic and to parabolic equations, see e.g.~\cite[pp. 20]{RainerPicard2013}. The general idea has been applied in a number of studies to integro-differential equations (\cite{Trostorff2012a}), fractional differential equations (\cite{R.Picard2012}), problems with impedance type boundary conditions (\cite{Picard2012}), electro-seismic waves (\cite{McGhee2011}) or differential-algebraic systems arising in control theory (\cite{Picard2013}).

The aim of this article is twofold. At first we show that the technique to prove well-posedness in \cite{RainerPicard2013} applies to a more general situation covering the solution theories given in \cite{RainerPicard2013,PicPhy,Picard2012}, which can then be understood as a unifying strategy to tackle well-posedness issues in linear problems of mathematical physics.

The solution theory is built up in a certain space of exponentially weighted functions. In applications to particular examples one might choose this exponential weight small enough. Note that this exponential weight can be thought of as a $L^2$-analogue of the exponential growth of solutions in semi-group theory. Since changing the weight goes along with different underlying Hilbert spaces the question arises, whether the solution theory depends on the weight. Thus, the second aim of the present article is to show -- roughly speaking -- that the solution theory barely depends on the weight provided the operators involved barely depend on the weight. The latter is precisely the theorem, which one would hope for. The latter issue has not been addressed in \cite{RainerPicard2013} but it applies to the situation mentioned there.

As this article is intended to substantiate the results previously found in \cite{RainerPicard2013}, we will not give concrete examples and only sketch possible applications in Example \ref{ex:TDIDE}. The main issue will be that the problems from the linear theory\footnote{We refer the reader to \cite{S.Trostorff,Trostorff2012,Trostorff2011} for possibilites to deal with maximal monotone (non-linear) relations.} well fit into the scheme developed. In order to do so, we provide some basic notions and definitions in Section \ref{sec:prel}. In this section some results from \cite{Picard,Waurick2013b,Waurick2013c} are summarized. Section \ref{sec:wpr} deals with the statement of the well-posedness theorem and elaborates the relations to the ones in \cite{Picard2012} and \cite{RainerPicard2013}. In Section \ref{sec:ind_of_nu} we show the aforementioned independence of the exponential weight and Section \ref{sec:proof_of_solth} provides the proof of the well-posedness theorem.

\section{Preliminaries}\label{sec:prel}

 For $\nu\in \R$ and a Hilbert space $H$, we denote by $L_\nu^2(\R;H)$ the space of (equivalence classes of) measurable $H$-valued functions with respect to the Lebesgue-weight $x\mapsto e^{-2\nu x}$, thus $L_0^2(\R;H)=L^2(\R;H)$. We use $\langle \cdot,\cdot\rangle_\nu$ and $\abs{\cdot}_\nu$ to denote the scalar product and the norm in $L_\nu^2(\R;H)$, respectively. If the value of $\nu$ is clear from the context, we drop the respective indices. Denoting by $H_{\nu,1}(\R;H)$ the space of $H$-valued weakly differentiable functions with weak derivative representable as a $L_\nu^2(\R;H)$-function, we can show that
\[
  \partial_{0,\nu} \colon H_{\nu,1}(\R;H) \subseteqq L_\nu^2(\R;H)\to L_\nu^2(\R;H),\phi\mapsto \phi'
\]
 defines a continuously invertible operator if $\nu\neq 0$. For $\nu>0$ and $f\in L_\nu^2(\R;H)$ the inverse of $\partial_{0,\nu}$ can be expressed with the help of the Bochner-integral
\[
  \partial_{0,\nu}^{-1}f(t) = \int_{-\infty}^t f(\tau)\dd \tau \quad (t\in\R).
\]
Moreover, the \emph{Fourier-Laplace transformation} $\s L_\nu$ being the (unitary) closure of the operator\footnote{The space of compactly supported continuous functions on $\R$ with values in $H$ is denoted by $C_c(\R;H)$.}
\begin{align*}
  C_c(\R;H) \subseteqq L_\nu^2(\R;H) &\to L^2(\R;H)\\
                                    f&\mapsto \left( x\mapsto \frac{1}{\sqrt{2\pi}}\int_{\R} e^{-ixy-\nu y}f(y)\dd y\right)
\end{align*}
is an explicit spectral representation for $\partial_{0,\nu}$. Indeed, denoting by $m\colon D(m)\subseteqq L^2(\R;H)\to L^2(\R;H)$ the multiplication-by-argument operator given by $mf(x)\coloneqq xf(x)$, $x\in \R$, for $f$ lying in the maximal domain $D(m)$ of $m$, one can prove the formula (\cite[p.~161-163]{Akhiezer1981})
\begin{equation}\label{eq:funcalc}
   \partial_{0,\nu} = \s L_\nu^* ( im+\nu) \s L_\nu.
\end{equation}
Again, if the value of $\nu$ is clear from the context, we drop the index in the notation of $\partial_{0,\nu}$.
 Equation \eqref{eq:funcalc} gives rise to a functional calculus for $\partial_{0,\nu}^{-1}$: If $M\colon B(r,r)\to L(H)$ is a bounded, analytic function defined on some complex ball with radius $r>\frac{1}{2\nu}$ centered at $r$ with values in the space of bounded linear operators on $H$, we define 
\[
   M\left(\partial_{0,\nu}^{-1}\right) \coloneqq \s L_\nu^* M\left(\frac{1}{im+\nu}\right)\s L_\nu,
\]
where $M\left(\frac{1}{im+\nu}\right)\phi(t)\coloneqq M\left(\frac{1}{it+\nu}\right)\phi(t)$ for $\phi\in L^2(\R;H)$, $t\in \R$. Invoking \cite[Theorem 6.5]{Thomas1997} (see also \cite{Jac2000,Par2004}), we realize that bounded and analytic functions of $\partial_{0,\nu}^{-1}$ in the aforementioned sense are precisely the ones being causal in the following sense: 
 \begin{Def}[causality, see e.g.~{\cite[Definition 3.1.47]{Picard}}] Let $H_0,H_1$ be Hilbert spaces, $\nu>0$, $M\colon D(M)\subseteqq L_\nu^2(\R;H_0)\to L_\nu^2(\R;H_1)$. We say that $M$ is \emph{causal} if for all $a\in \R$ and $f,g\in D(M)$ with\footnote{For a function $\psi\colon \R \to \mathbb{C}$ we denote by $\psi(m_0)$ the associated multiplication operator acting on Hilbert space-valued functions $f\colon \R\to H$ in the way that $\left(\psi(m_0)f\right)(t)\coloneqq \psi(t)f(t)$, $t\in \R$.} $\1_{\R_{\leqq a}}(m_0)(f-g)=0$ then $\1_{\R_{\leqq a}}(m_0)(M(f)-M(g))=0$. 
 \end{Def}
The next lemma is almost immediate:
\begin{Le}\label{le:operator_equality_causality} Let $H$ be a Hilbert space, $M\colon D(M)\subseteqq L_\nu^2(\R;H)\to L_\nu^2(\R;H)$. Assume that $\1_{\R_{\leqq a}}(m_0)[D(M)]\subseteqq D(M)$ for all $a\in\R$. Then the following statements are equivalent:
\begin{enumerate}[(i)]
 \item\label{char_caus1} $M$ is causal;
 \item\label{char_caus2} for all $a\in \R$ we have $\1_{\R_{\leqq a}}(m_0)M = \1_{\R_{\leqq a}}(m_0)M\1_{\R_{\leqq a}}(m_0)$.
\end{enumerate} 
\end{Le}
In \cite{Waurick2013b} it has been found that at least for closure procedures the latter concept of causality seems not to be appropriate. A possible way to overcome this difficulty is to introduce the following notion:
\begin{Def}[strong causality, {\cite{Waurick2013b}}] Let $H_0,H_1$ be Hilbert spaces, $\nu>0$, $M\colon D(M)\subseteqq L_\nu^2(\R;H_0)\to L_\nu^2(\R;H_1)$. We say that $M$ is \emph{strongly causal} if for all $R>0$, $a\in \R$, $\phi\in L_\nu^2(\R;H_1)$ the mapping \begin{align*}
               \left( B_{M}(0,R), \abs{ \1_{\R_{\leqq a}}(m_0)\left(\cdot-\cdot\right)}\right) &\to \left(L_\nu^2(\R;H), \abs{\langle \1_{\R_{\leqq a}}(m_0)\left(\cdot-\cdot\right),\phi\rangle}\right) \\
   f&\mapsto Mf,
\end{align*} 
is uniformly continuous, where $B_M(0,R)\coloneqq \{ f\in D(M); \abs{f}^2+\abs{Mf}^2< R^2\}$.
\end{Def}
The latter notion is equivalent to causality for particular situations.
\begin{Sa}[{\cite[Theorem 1.6]{Waurick2013b}}]\label{thm:causal_cont_causal} Let $H_0, H_1$ be Hilbert spaces, $\nu>0$, $M\colon D(M)\subseteqq L_\nu^2(\R;H_0)\to L_\nu^2(\R;H_1)$ densely defined, linear and closable. Then the following statements are equivalent:
\begin{enumerate}[(i)]
 \item\label{caus_cont_caus1} $\overline{M}$ is causal;
 \item\label{caus_cont_caus2} $M$ is strongly causal.
\end{enumerate}
\end{Sa}
In \cite{RainerPicard2013}, we found that multiplication operators, which -- in contrast to functions of $\partial_{0,\nu}^{-1}$ -- are not time-translation invariant, can also be dealt with as coefficients in a solution theory for certain linear evolutionary equations. A notion containing both the aforementioned functions of $\partial_{0,\nu}^{-1}$ as well as multiplication operators is the one given in \cite{Waurick2013c}, the notion of evolutionary mappings:
\begin{Def}[evolutionary mappings, {\cite[Definition 2.1]{Waurick2013c}}]
Let $H_{0},H_{1}$ be Hilbert spaces, $\nu_1>0.$ We call a linear
mapping 
\begin{equation}
M\colon D(M)\subseteqq\bigcap_{\nu\geqq\nu_1}L_{\nu}^{2}(\mathbb{R};H_{0})\to\bigcap_{\nu\geqq\nu_{1}}L_{\nu}^{2}(\mathbb{R};H_{1})\label{eq:evolutionary_map}
\end{equation}
\emph{evolutionary (at $\nu_{1}$)} 
if $D(M)$ is dense in $L_{\nu}^{2}(\mathbb{R};H_0)$  and $M$ is a closable operator from $L_\nu^2(\R;H_0)$ to $L_\nu^2(\R;H_1)$ for all $\nu\geqq\nu_{1}$.
$M$ is called \emph{bounded}, if, in addition, $M$ extends to a bounded linear operator from $L_{\nu}^{2}(\mathbb{R};H_{0})$
to $L_{\nu}^{2}(\mathbb{R};H_{1})$ for all $\nu\geqq\nu_{1}$ such that%
\footnote{For a linear operator $A$ from $L_{\nu}^{2}(\mathbb{R};H_{0})$ to
$L_{\nu}^{2}(\mathbb{R};H_{1})$ we denote its operator norm by $\left\Vert A\right\Vert _{L(L_{\nu}^{2}(\mathbb{R};H_{0}),L_{\nu}^{2}(\mathbb{R};H_{1}))}$.
If the spaces $H_{0}$ and $H_{1}$ are clear from the context, we
shortly write $\left\Vert A\right\Vert _{L(L_{\nu}^{2})}$. %
} 
\[
\limsup_{\nu\to\infty}\left\Vert M\right\Vert _{L(L_{\nu}^{2}(\mathbb{R};H_{0}),L_{\nu}^{2}(\mathbb{R};H_{1}))}<\infty.
\]
 For evolutionary mappings $M$, the closure of $M$ in some $L_{\nu}^{2}$ will
be denoted by $M_\nu$.
We define the sets
\[
\s C_{\text{ev},\nu_{1}}(H_{0},H_{1})\coloneqq\{M;M\text{ is as in (\ref{eq:evolutionary_map}) and is evolutionary at }\nu_{1}\}
\]
and
\[
  L_{\text{ev},\nu_{1}}(H_{0},H_{1})\coloneqq\{M;M\text{ is as in (\ref{eq:evolutionary_map}), is evolutionary at }\nu_{1}\text{ and bounded}\}
\]
We abbreviate $\s C_{\text{ev},\nu_{1}}(H_{0})\coloneqq \s C_{\text{ev},\nu_{1}}(H_{0},H_{0})$ and $L_{\text{ev},\nu_{1}}(H_{0})\coloneqq L_{\text{ev},\nu_{1}}(H_{0},H_{0})$.
\end{Def}

\begin{Bei}\label{ex:bddfunctionsareevolutionary} Let $H$ be a Hilbert space, $r>0$, $M\colon B(r,r)\to L(H)$ bounded and analytic. For $\nu>\frac{1}{2r}$ the operator $M\left(\partial_{0}^{-1}\right)$ is a bounded evolutionary mapping at $\nu$. Indeed, take $D\left(M\left(\partial_{0}^{-1}\right)\right)=C_{\infty,c}(\R;H)$, the space of compactly supported $H$-valued functions, which are indefinitely differentiable. Then $C_{\infty,c}(\R;H)$ is dense in $L_{\nu'}^2(\R;H)$ for every $\nu'\geqq \nu$. Moreover, Cauchy's integral theorem shows that 
\[
   \s L_{\nu'}^* M\left(\frac{1}{im+\nu'}\right)\s L_{\nu'}\phi = \s L_{\nu}^* M\left(\frac{1}{im+\nu}\right)\s L_{\nu}\phi 
\]
for every $\nu'\geqq \nu$ and $\phi\in C_{\infty,c}(\R;H)$, see also \cite[Lemma 3.6]{Trostorff2013a} for more details. Moreover, the operator norm of $M\left(\partial_0^{-1}\right)$ in $L_{\nu'}^2(\R;H)$ is easily estimated by $\sup_{z\in B(r,r)} \Abs{M(z)}_{L(H)}$.
\end{Bei}

\begin{Bei}[{\cite[Example 2.3]{Waurick2013c}}]\label{ex:operator-valued-multiplication}
 Let $H$ be a Hilbert space and let $L_{s}^{\infty}(\mathbb{R};L(H))$
be the space of bounded strongly measurable functions from $\mathbb{R}$
to $L(H)$. For $A\in L_{s}^{\infty}(\mathbb{R};L(H))$ we denote
the associated multiplication operator on $L_{\nu}^{2}(\mathbb{R};H)$
by $A(m_{0})$. It is easily verified that $A(m_{0})\in\bigcap_{\nu>0}L_{\text{ev},\nu}(H)$.
\end{Bei}

\begin{rem}\label{evo=caus} Evolutionarity and causality are rather closely connected. Indeed, let $M\in  L_{\text{ev},\nu_{1}}(H_{0},H_{1})$ and assume that for all $a\in \R$ the set $D(M\1_{\R_{\leqq a}}(m_0))\cap D(M)$ is dense in $L_\nu^2(\R;H_0)$ for $\nu\geqq\nu_1$\footnote{This assumption is met if, for instance, $D(M)$ contains $C_{\infty,c}(\R;D)$ for $D$ being a Hilbert space densely embedded into $H_0$. An example, where the intersection $D(M\1_{\R_{\leqq a}}(m_0))\cap D(M)$ only contains $0$ for all $a\in \R$ is given in \cite[Example 1.4]{Waurick2013b}.}. Then $M_\nu$ is causal for all $\nu\geqq\nu_1$. Indeed, let $f\in L_\nu^2(\R;H_0)$, $a\in\R$ and assume that $\1_{\R_{\leqq a}}(m_0)f=0$. We choose $(\phi_n)_{n}$ in $D(M\1_{\R_{\leqq a}}(m_0))\cap D(M)$ such that $\phi_n\to f$ in $L_\nu^2(\R;H_0)$ as $n\to\infty$. For all $n\in\N$ we have $\1_{\R_{\geqq a}}(m_0)\phi_n=\phi_n-\1_{\R_{< a}}(m_0)\phi_n \in D(M)$ and 
\[
   \1_{\R_{\geqq a}}(m_0)\phi_n \to \1_{\R_{\geqq a}}(m_0)f= f\quad (n\to\infty)
\]
 in $L_\nu^2(\R;H_0)$. In particular, the latter implies that $\psi_n\coloneqq \1_{\R_{\geqq a}}(m_0)\phi_n$ approximates $f$ in $L_{\nu'}^2(\R;H_0)$ for all $\nu'\geqq \nu$ and, thus, $\psi_n\to f$ in $L_{\textnormal{loc}}^2(\R;H)$ as $n\to\infty$. Now, we follow the idea of \cite[Proof of Theorem 4.5]{Kalauch}. For this let $\phi\in C_{\infty,c}(\R;H_1)$ with support bounded above by $a$. For $\nu'\geqq\nu$ we get that 
\begin{align*}
   \abs{\langle M_\nu f , \phi \rangle_{0}} & = \lim_{n\to\infty}\abs{\langle M_\nu \psi_n,\phi\rangle } \\
                                            & = \lim_{n\to\infty}\abs{\langle M \psi_n,\phi\rangle } \\
                                            & \leqq \lim_{n\to\infty}\abs{M \psi_n}_{\nu'}\abs{\phi}_{-\nu'} \\
                                            & \leqq \Abs{M}_{L(L_{\nu'}^2)} \abs{f}_{\nu'}\abs{\phi}_{0}e^{\nu'a} \\
					    & = \Abs{M}_{L(L_{\nu'}^2)} \abs{f(\cdot +a)}_{\nu'}\abs{\phi}_{0}  
\end{align*}
Letting $\nu'$ tend to infinity, we deduce that $\langle M_\nu f , \phi \rangle_{0}=0$. Hence, $M_\nu f=0$ on $(-\infty,a]$. 
\end{rem}

\section{The well-posedness result}\label{sec:wpr}

The well-posedness result will be formulated within the following situation. Throughout, let $H$ be a Hilbert space and $\nu>0$.

\begin{hyp}[on the material law operator]\label{hyp:mat_law} Let $\s M, \s N \in L(L_\nu^2(\R;H))$. Assume that there exists $M\in L(L_\nu^2(\R;H))$ such that 
\[
   \s M \partial_0 \subseteqq \partial_0 \s M + M.
\] 
Further assume that both $\s M$ and $\s N$ are causal.
\end{hyp}

\begin{hyp}[on the unbounded spatial operator]\label{hyp:A} Let $\s A\colon D(\s A)\subseteqq L_\nu^2(\R;H)\to L_\nu^2(\R;H)$ be densely defined, closed, linear and such that $\partial_0^{-1}\s A \subseteqq \s A \partial_0^{-1}$.
\end{hyp}

\begin{rems}
 (a) Note that Hypothesis \ref{hyp:A} reflects the fact that we will treat autonomous $\s A$. Since we do not assume this commutativity condition for the material law operator, we allow for operators, which are not time-translation invariant in Hypothesis \ref{hyp:mat_law}.

 (b) The assumption in Hypothesis \ref{hyp:A} implies that $(1+\eps\partial_0)^{-1}\s A \subseteqq \s A (1+\eps\partial_0)^{-1}$ for every $\eps>0$, since $(1+\eps\partial_0)^{-1}$ is a bounded (continuous) Borel function of $\partial_0^{-1}$.

 (c) Hypothesis \ref{hyp:A} particularly implies that $D(\partial_0)\cap D(\s A)$ is dense in $L_\nu^2(\R;H)$. Indeed, let $\eps>0$ and let $f\in D(\s A)$. Then, by (b), we have $(1+\eps\partial_0)^{-1}\s A \subseteqq\s A(1+\eps\partial_0)^{-1}$. Thus, $(1+\eps\partial_0)^{-1}f \in D(\partial_0)\cap D(\s A)$. Since $(1+\eps\partial_0)^{-1}f\to f$ as $\eps\to 0$ (see e.g.~\cite[Lemma 2.7]{RainerPicard2013}), we conclude that $D(\partial_0)\cap D(\s A)$ is dense in $D(\s A)$. The density of $D(\s A)$ in $L_\nu^2(\R;H)$ yields the assertion.
\end{rems}

The results reads as follows.

\begin{Sa}\label{thm:genSolth} Let $\s M,\s N,\s A$ be as in Hypotheses \ref{hyp:mat_law} and \ref{hyp:A}. Assume there exists $c>0$ such that the positivity conditions
\begin{equation}\label{ineq:pos_cutoff}
 \Re\langle \left(\partial_0\s M+\s N +\s A\right)\phi,\1_{\R_{\leqq a}}(m_0)\phi\rangle \geqq c \langle \phi,\1_{\R_{\leqq a}} (m_0)\phi\rangle
\end{equation}
and 
\begin{equation}\label{ineq:pos_adjoint}
   \Re\langle \left(\left(\partial_0\s M+\s N\right)^* +\s A^*\right)\psi,\psi\rangle \geqq c \langle \psi,\psi\rangle
\end{equation}
  hold for all $a\in \R$, $\phi\in D(\partial_0)\cap D(\s A)$, $\psi\in D(\partial_0)\cap D(\s A^*)$.

Then $\s B\coloneqq \overline{\partial_0\s M+\s N+\s A}$ is continuously invertible, $\Abs{\s B^{-1}}\leqq \frac{1}{c}$, and the operator $\s B^{-1}$ is causal in $L_\nu^2(\R;H)$.
\end{Sa}

For later purposes, we also have the following density result:

\begin{Le}\label{le:nice_density} Under the hypothesis of Theorem \ref{thm:genSolth}, let $F\subseteqq D(\s A)$ be a core for $\s A$. Then $\s F\coloneqq \lin\bigcup_{\delta>0}\left(1+\delta\partial_0\right)^{-1}[F]$ is a core for $\overline{{\partial_0 \s M+\s N+\s A}}$
\end{Le}

\begin{rems} 

(a) Of course there is also an adapted perturbation theorem, similar to the ones given in \cite[Theorem 2.17 and Theorem 2.19]{RainerPicard2013}. The proof of the respective results will be obvious from the proofs in \cite{RainerPicard2013}. Thus, we will not repeat them here. 

(b) The condition on $\s A$ to commute with $\partial_0^{-1}$ can be relaxed in the sense that it suffices to assume that $\s A$ and $\partial_0$ have a bounded commutator. We will address this non-commutativity relation in a future article.

(c) The truncation in the positive definiteness condition \eqref{ineq:pos_cutoff} is needed in order to obtain causality. The proof in Section \ref{sec:proof_of_solth} will show that the well-posedness result, i.e., continuous invertibility of the closure of $\partial_0\s M+\s N+\s A$, can also be obtained if one assumes  
\begin{equation}\label{ineq:pos_cutoff2}
 \Re\langle \left(\partial_0\s M+\s N +\s A\right)\phi,\phi\rangle \geqq c \langle \phi,\phi\rangle
\end{equation}
instead of \eqref{ineq:pos_cutoff} for all $\phi\in D(\partial_0)\cap D(\s A)$.
\end{rems}

 We will postpone a proof of Theorem \ref{thm:genSolth} to Section \ref{sec:proof_of_solth}. In the remainder of this section, we show that Theorem \ref{thm:genSolth} contains \cite[Theorem 2.13]{RainerPicard2013} and the solution theory stated in \cite{Picard2012} and \cite[Theorem 6.2.5]{Picard} as special cases. In order to obtain the latter, we observe the following consequence of Theorem \ref{thm:genSolth}:

\begin{Sa}\label{thm:Solth}
 Let $\s M,\s N,\s A$ be as in Hypotheses \ref{hyp:mat_law} and \ref{hyp:A}. Assume in addition that $\s A$ is maximal monotone in $L_\nu^2(\R;H)$, $\s A$ is causal and $\1_{\R_{\leqq 0}}(m_0)[D(\s A)]\subseteqq D(\s A)$.
Moreover, assume the positive definiteness condition
\begin{equation}\label{ineq:thm:Solth}
   \Re \left\langle \left(\partial_0 \s M+\s N\right)u,\1_{\R_{\leqq a}}(m_0)u\right\rangle_{\nu}\geqq c\left\langle u,\1_{\R_{\leqq a}}(m_0)u \right\rangle_{\nu} 
\end{equation}
for all $u\in D(\partial_0)$, $a\in \R$ and some $c>0$.
Then $0\in \rho\left(\overline{\partial_0 \s M+\s N+\s A}\right)$ and the operator $\left(\overline{\partial_0 \s M+\s N+\s A}\right)^{-1}$ is causal.
\end{Sa}

Next we show that Theorem \ref{thm:Solth} follows from Theorem \ref{thm:genSolth}. Indeed, this follows from the next two lemmas:

\begin{Le}\label{le:ineq_cutoff} Under the hypotheses of Theorem \ref{thm:Solth}, we get for every $a\in \R$ and $\phi\in D(\partial_0)\cap D(\s A)$ that
\[
  \Re \langle \left(\partial_0\s M+\s N + \s A\right)\phi,\1_{\R_{\leqq a}}(m_0)\phi\rangle\geqq c\langle\phi,\1_{\R_{\leqq a}}(m_0)\phi\rangle.
\]
\end{Le}
 
\begin{Le}\label{le:ineq_adjoint} Under the hypotheses of Theorem \ref{thm:Solth}, we have for all $\psi \in D(\partial_0)\cap D(\s A)$
\[
   \Re \langle \left( \left(\partial_0 \s M + \s N\right)^* + \s A^* \right)\psi,\psi\rangle \geqq c\langle\psi,\psi\rangle.
\]
\end{Le}

The proofs of the Lemmas \ref{le:ineq_adjoint} and \ref{le:ineq_cutoff} need the following preliminary observation:
\begin{Le}\label{lem:positiv_of_adjoint} Under the hypotheses of Theorem \ref{thm:Solth}, we have for all $\phi\in D(\partial_0)$:
\[
   \Re \langle \left(\partial_0\s M +\s N\right)\phi,\phi\rangle = \Re \langle \left(\partial_0\s M +\s N\right)^*\phi,\phi\rangle \geqq c\langle\phi,\phi\rangle.
\]
\end{Le}
\begin{proof}
 We observe that the equality is true for all $\phi \in D(\partial_0\s M+\s N)\cap D(\left(\partial_0 \s M+\s N\right)^*)\cap D(\partial_0)$. The inclusion $D(\partial_0\s M+\s N)=D(\partial_0\s M)\supseteqq D(\partial_0)$ being obvious, the only thing left to prove is $D(\left(\partial_0 \s M+\s N\right)^*)\supseteqq D(\partial_0)$. For this, we compute (use $\partial_{0}^*=-\partial_0+2\nu$) that 
\[
   \left(\partial_0 \s M+\s N\right)^*=\left(\partial_0 \s M\right)^*+\s N^* \supseteqq \s M^*\partial_0^*+\s N^* = -\s M^*\partial_0 + 2\nu \s M^* + \s N^*, 
\]
 which yields the assertion.
\end{proof}

\begin{proof}[Proof of Lemma \ref{le:ineq_cutoff}]
  Let $a\in\R$ and $\phi\in D(\s A)\cap D(\partial_0)$. We compute
\begin{align*}
  & \Re \langle \left(\partial_0 \s M+\s N+\s A\right)\phi,\1_{\R_{\leqq a}}(m_0)\phi\rangle \\
  & =\Re \langle \left(\partial_0 \s M+\s N\right)\phi+\1_{\R_{\leqq a}}(m_0)\s A\phi,\1_{\R_{\leqq a}}(m_0)\phi\rangle\\
  & =\Re \langle \left(\partial_0 \s M+\s N\right)\phi+\1_{\R_{\leqq a}}(m_0)\s A\1_{\R_{\leqq a}}(m_0)\phi,\1_{\R_{\leqq a}}(m_0)\phi\rangle\\
  & \geqq c \langle \phi,\1_{\R_{\leqq a}}(m_0)\phi\rangle+\Re\langle \s A\1_{\R_{\leqq a}}(m_0)\phi,\1_{\R_{\leqq a}}(m_0)\phi\rangle\\
  & \geqq c \langle \phi,\1_{\R_{\leqq a}}(m_0)\phi\rangle,
\end{align*}
where we have used the monotonicity of $\s A$ and that $\1_{\R_{\leqq a}}(m_0)$ leaves the domain of $\s A$ invariant. Indeed, the inclusion $\1_{\R_{\leqq a}}(m_0)[D(\s A)]\subseteqq D(\s A)$ for all $a\in \R$ follows from $\1_{\R_{\leqq 0}}(m_0)[D(\s A)]\subseteqq D(\s A)$. For this, let $a\in\R$ and denote by $\tau_a\in L(L_\nu^2(\R;H))$ the shift of a function $f$ by $a$, i.e., $\tau_af\coloneqq f(\cdot+a )$. Since $\tau_a$ is a Borel function of $\partial_0^{-1}$ it thus commutes with $\s A$. In particular, we have for all $\phi\in D(\s A)$ also $\tau_a \phi \in D(\s A)$. Hence,  $\1_{\R_{\leqq a}}(m_0) \phi = \1_{\R_{\leqq a}}(m_0)\tau_{-a}\tau_a \phi = \tau_{-a}\1_{\R_{\leqq 0}}(m_0)\tau_a \phi \in D(\s A)$, since $\1_{\R_{\leqq 0}}(m_0)\tau_a \phi\in D(\s A)$, by assumption and $\tau_{-a}[D(\s A)]\subseteqq D(\s A)$. 
\end{proof}

\begin{proof}[Proof of Lemma \ref{le:ineq_adjoint}]
  It is well-known that the maximal monotonicity of $\s A$ implies the same for $\s A^*$, see e.g.~\cite[Exercise 1.7.5]{Arendt2005/2006}. (Use that maximal monotonicity for $\s A$ is equivalent to $(0,\infty)\subseteqq \rho(-\s A)$ and $\Abs{\lambda (\lambda+A)}\leqq 1$ for all $\lambda\in\R_{>0}$.) Hence, the result follows from Lemma \ref{lem:positiv_of_adjoint}.
\end{proof}

Now, we discuss the relationship of Theorem \ref{thm:genSolth} to the well-posedness theorems stated in \cite[Solution Theory]{PicPhy} (\cite[Theorem 6.2.5]{Picard}) with its generalization in  \cite[Theorem 1.2]{Picard2012} and its time-dependent analogues \cite[Theorem 2.13 and Theorem 2.15]{RainerPicard2013}.

\subsection*{Relationship to the well-posedness condition in \cite[Theorem 1.2]{Picard2012}}

The well-posedness theorem in \cite{Picard2012} has been applied to a model for acoustic waves with impedance type boundary conditions. The theorem reads as follws.

\begin{Sa}[{\cite[Theorem 1.2]{Picard2012}}] \label{thm:well_iwota} Let $H$ be a Hilbert space, $\nu>0$, $r>\frac{1}{2\nu}$. Denote by $M\colon B(r,r)\to L(H)$ a bounded and analytic function and let $\s A : D(\s A)\subseteqq L_\nu^2(\R;H)\to L_\nu^2(\R;H)$ satisfy Hypothesis \ref{hyp:A}. Assume that
 \[
    \Re \langle \left(\partial_0 M\left(\partial_0^{-1}\right)+\s A\right)\phi,\1_{\R_{\leqq 0}}(m_0)\phi\rangle \geqq c\langle \phi,\1_{\R_{\leqq 0}}(m_0)\phi\rangle,
 \]as well as
 \[
    \Re \langle \left(\left(\partial_0 M\left(\partial_0^{-1}\right)\right)^*+\s A^*\right)\psi,\psi\rangle \geqq c\langle \psi,\psi\rangle.
 \]
for all $\phi\in D(\partial_0)\cap D(\s A)$ and $\psi\in D(\partial_0)\cap D(\s A^*)$.
  Then $\overline{\partial_0M(\partial_0^{-1})+\s A}$ is continuously invertible and causal in $L_\nu^2(\R;H)$.
\end{Sa}

Regarding the positive definiteness conditions in Theorem \ref{thm:well_iwota}, we realize that these conditions are precisely the ones in Theorem \ref{thm:genSolth} (with $\s M=M(\partial_0^{-1})$, $\s N=0$). Indeed, since both $M\left(\partial_0^{-1}\right)$ and $\s A$ commute with time-translation, the positive definiteness condition in Theorem \ref{thm:well_iwota} carries over to the situation, where $\1_{\R_{\leqq 0}}(m_0)$ is replaced by $\1_{\R_{\leqq a}}(m_0)$, $a\in \R$. For the latter also see the argument in the proof of Lemma \ref{le:ineq_cutoff}. The only thing left to prove is that  $M\left(\partial_0^{-1}\right)$ has a bounded commutator with $\partial_0$, which is obvious since the former is a function of the latter. Indeed, in this situation the commutator equals $0$.

\subsection*{Relationship to the well-posedness condition in \cite[Theorem 2.13 and Theorem 2.15]{RainerPicard2013}}

Recall the well-posedness result in \cite{RainerPicard2013}:

\begin{Sa}[{\cite[Theorem 2.15]{RainerPicard2013}}]\label{thm:maria} Let $H$ be a Hilbert space, $\nu>0$, $A\colon D(A)\subseteqq H\to H$ skew-selfadjoint\footnote{We use the canonical extension of $A$ as multiplication operator in $L_\nu^2(\R;H)$ and use the same notation.}. Let $M_0,M_1\in L_s^\infty(\R;L(H))$ (see Example \ref{ex:operator-valued-multiplication} for a definition). 
Assume in addition that $M_0$ is strongly differentiable, Lipschitz continuous and that\footnote{For a bounded linear operator $B\in L(H)$, we define its (selfadjoint) real-part by $\Re B \coloneqq \frac{1}{2}(B+B^*)$. Then $\Re B \geqq c$ means that $\Re B-c$ lies in the cone $\s K\subseteqq L(H)$ of positive definite bounded linear operators, i.e., $T\in \s K$ if $\langle T\phi,\phi\rangle\geqq 0$ for all $\phi\in H$.} 
\begin{equation}\label{ineq:pos_def_time}
    \nu M_0(t)+ \frac{1}2 \dot M_0(t)+\Re M_1(t) \geqq c
\end{equation}
in $L(H)$ for almost every $t\in \R$ and all sufficiently large $\nu$.
Then $\overline{\partial_0M_0(m_0)+M_1(m_0)+ A}$ is continuously invertible and causal in $L_\nu^2(\R;H)$ (if $\nu$ is sufficiently large).
 \end{Sa}

For this theorem it should be noted that this is a direct consequence of Theorem \ref{thm:Solth} (applied to $M_0(m_0)=\s M$, $M_1(m_0)=\s N$ and $A=\s A$). Indeed, the strong differentiability and the Lipschitz continuity of $M_0$ ensure that $M_0(m_0)$ and $\partial_0$ have a bounded commutator, see \cite[Lemma 2.1]{RainerPicard2013}. Moreover, skew-selfadjoint operators are maximal monotone and condition \eqref{ineq:pos_def_time} implies condition \eqref{ineq:thm:Solth} assumed in Theorem \ref{thm:Solth}, see \cite[Lemma 2.6]{RainerPicard2013}. The other conditions are easily verified. Indeed, using Remark \ref{evo=caus} and taking into account that multiplication operators induced by bounded strongly measurable mappings are evolutionary at every $\nu\geqq 0$ (see Example \ref{ex:operator-valued-multiplication}) and defined on the whole of $L_\nu^2(\R;H)$ for every $\nu\geqq 0$, we realize that these multiplication operators are strongly causal. Note that a similar argument also applies to $A$, which is a bounded evolutionary mapping in $L_{\textnormal{ev},\nu}(D_A,H)$ for every $\nu\geqq 0$, where $D_A$ denotes the domain of $A$ endowed with the graph norm of $A$.

\subsection*{An example of particular material laws}

In the introduction we elaborated the applicability of the theorems above and cited the respective references. Thus, there is no need to repeat them here. Due to the two latter observations concerning the entailment of the well-posedness theorems given in the literature, we give an example of material laws, which are neither covered by the one or the other well-posed theorem but by the 
well-posedness theorem discussed in this article.

\begin{Bei}[Time-dependent integro-differential-algebraic equations]\label{ex:TDIDE} Let $\nu>0$, $r>\frac{1}{2\nu}$.  Assume $H=H_0\oplus H_1$ and take $M\colon B(r,r)\to L(H_0)$ analytic and bounded. Moreover, let $M_0,M_1 \in L_s^\infty(\R;L(H_1))$ satisfy the conditions in Theorem \ref{thm:maria} with $H$ replaced by $H_1$ and assume that $\Re \partial_0M(\partial_0^{-1})\geqq c$ in $L_\nu^2(\R;L(H_0))$. Then, for \emph{any} maximal monotone operator $A\colon D(A)\subseteqq H\to H$, the operator sum
\[
   \partial_0 \begin{pmatrix}
    M(\partial_0^{-1}) & 0 \\ 0 & M_0(m_0) 
   \end{pmatrix} + \begin{pmatrix} 0 & 0 \\ 0& M_1(m_0)\end{pmatrix} + A
\]
is closable with continuous invertible closure for sufficiently large $\nu$. Since the block structure of $\partial_0 \begin{pmatrix}
    M(\partial_0^{-1}) & 0 \\ 0 & M_0(m_0) 
   \end{pmatrix} + \begin{pmatrix} 0 & 0 \\ 0& M_1(m_0)\end{pmatrix}$ and $A$ need not be comparable, the continuous invertibility does \emph{not} follow either from Theorem \ref{thm:well_iwota} or \ref{thm:maria} or both of them. Hence, we need to invoke Theorem \ref{thm:genSolth}. 
\end{Bei}

%

\section{On the independence of $\nu$}\label{sec:ind_of_nu}

Of course, a natural question in the general setting of Theorem \ref{thm:genSolth} and more particularly in Theorem \ref{thm:maria} is whether the solution operator depends an the particular choice of $\nu$. In \cite[Theorem 2.13]{RainerPicard2013} the independence of $\nu$ has not been adressed. The next theorem shows that the solution indeed does not depend on the parameter $\nu$ in the following sense:

\begin{Sa}\label{thm:independence_of_nu}
 Let $H$ be a Hilbert space, $\nu_0>0$, $\s M,\s N\in L_{\textnormal{ev},\nu_0}(H)$, $\s A \in \s C_{\textnormal{ev},\nu_0}(H)$. 
 Assume that $\s M_\nu,\s N_\nu,\s A_\nu$ satisfy Hyptheses \ref{hyp:mat_law} and \ref{hyp:A}, respectively, for all $\nu\geqq \nu_0$. Moreover, assume that\[
   \Re\langle \left(\partial_{0,\nu}\s M_\nu+\s N_\nu +\s A_\nu\right)\phi,\1_{\R_{\leqq a}}(m_0)\phi\rangle_\nu \geqq c \langle \phi,\1_{\R_{\leqq a}} (m_0)\phi\rangle_\nu 
\]
and 
\[
   \Re\langle \left(\left(\partial_{0,\nu}\s M_\nu+\s N_\nu\right)^* +\s A_\nu^*\right)\psi,\psi\rangle_\nu \geqq c \langle \psi,\psi\rangle_\nu 
\]  hold for all $a\in \R$, $\phi\in D(\partial_{0,\nu})\cap D(\s A_\nu)$, $\psi\in D(\partial_{0,\nu})\cap D(\s A_\nu^*)$, $\nu\geqq \nu_0$.
 Then for $S_\nu \coloneqq \left(\overline{\partial_{0,\nu} \s M_\nu+\s N_\nu+\s A_\nu}\right)^{-1}$ we have $\Abs{S_\nu}\leqq \frac{1}{c}$ and 
%
for all $\nu_1,\nu\geqq \nu_0$ that
\[
   S_\nu f = S_{\nu_1} f \quad\left(f\in L_{\nu_1}^2(\R;H)\cap L_\nu^2(\R;H)\right). 
\]
\end{Sa}



To begin with, we observe the following relationship of causality and the independence of $\nu$:

\begin{Le}\label{le:pre_indep} Let $\nu_1\geqq \nu_0\geqq 0$, $H_0, H_1$ Hilbert spaces, $M_i \in L(L_{\nu_i}^2(\R;H_0),L_{\nu_i}^2(\R;H_0))$ causal, $i\in\{0,1\}$. Assume that there exists $D\subseteqq L_{\nu_0}^2(\R;H_0)\cap L_{\nu_1}^2(\R;H_0)$ dense in $L_{\nu_0}^2(\R;H_0)$. Then $M_0$ and $M_1$ coincide on $L_{\nu_0}^2(\R;H_0)\cap L_{\nu_1}^2(\R;H_0)$.
\end{Le}
\begin{proof}
   Let $f\in L_{\nu_1}^2(\R;H_0)\cap L_\nu^2(\R;H_0)$. Let $a \in \R$. By definition, there exists $(\phi_n)_n$ in $D$ such that $\phi_n\to f$ in $L_{\nu_0}^2(\R;H_0)$ as $n\to\infty$. Moreover, we get that $\1_{\R_{\leqq a}}(m_0)\phi_n \to \1_{\R_{\leqq a}}(m_0)f$ in $L_{\nu_1}^2(\R;H_0)$ as $n\to\infty$ 
 and that $\left(M_{1}\left(\1_{\R_{\leqq a}}(m_0)(\phi_n - f)\right)\right)_n$ tends to $0$ in $L_{\nu_1}^2(\R;H_1)$.  Using that both $M_0$ and $M_{1}$ are everywhere defined and Lemma \ref{le:operator_equality_causality}, we deduce for $n\in\N$ and $\psi\in C_{\infty,c}(\R;H)$ that
 \begin{align*}
   & \abs{\langle \1_{\R_{\leqq a}}(m_0) \left(M_{0} f-M_{1}f\right),\psi\rangle_{\nu_0}} \\
   &\leqq \abs{\langle \1_{\R_{\leqq a}}(m_0) \left(M_{0} f-M_{0}\phi_n\right),\psi\rangle_{\nu_0}} +\abs{\langle \1_{\R_{\leqq a}}(m_0) \left(M_{1} f-M_{1}\phi_n\right),\psi\rangle_{\nu_1}} \\
   &\leqq \abs{\langle \1_{\R_{\leqq a}}(m_0) \left(M_{0} f-M_{0}\phi_n\right),\psi\rangle_{\nu_0}} \\ &\quad +\abs{\langle \1_{\R_{\leqq a}}(m_0) \left(M_{1} \1_{\R_{\leqq a}}(m_0)f-M_{1}\1_{\R_{\leqq a}}(m_0)\phi_n\right),\psi\rangle_{\nu_0}} \\
   &=\abs{\langle \1_{\R_{\leqq a}}(m_0) \left(M_{0} f-M_{0}\phi_n\right),\psi\rangle_{\nu_0}} \\ &\quad +\abs{\langle \1_{\R_{\leqq a}}(m_0) \left(M_{1} \1_{\R_{\leqq a}}(m_0)f-M_{1}\1_{\R_{\leqq a}}(m_0)\phi_n\right),e^{2\left(\nu_1-\nu_0\right)(\cdot)}\psi\rangle_{\nu_1}}\\
   &\to 0 \quad (n\to\infty).
\end{align*}
Thus, $M_0 f = M_{1}f$ almost everywhere. 
\end{proof}

We observe the following consequence of Lemma \ref{le:pre_indep}, which gives more insight to bounded (causal) evolutionary mappings:

\begin{Le}\label{le:indep_of_nu_for_bounded_mappings} Let $\nu\geqq\nu_1\geqq \nu_0$, $H_0,H_1$ Hilbert spaces. Let $M\in L_{\textnormal{ev},\nu_0}(H_0,H_1)$ be strongly causal. Then $M_\nu$ and $M_{\nu_1}$ coincide on $L_{\nu_1}^2(\R;H_0)\cap L_\nu^2(\R;H_0)$. 
\end{Le}
\begin{proof}
  The assertion follows from Theorem \ref{thm:causal_cont_causal} and Lemma \ref{le:pre_indep}.
\end{proof}

We come to the proof of Theorem \ref{thm:independence_of_nu}.

\begin{proof}[Proof of Theorem \ref{thm:independence_of_nu}] For $\nu\geqq \nu_0$ we define $\tilde S_\nu \coloneqq \left(\partial_{0,\nu}\s M_\nu+\s N_\nu+\s A_\nu\right)^{-1}$. From Theorem \ref{thm:genSolth}, we deduce that $\Abs{\tilde S_\nu}_{L(L_\nu^2)}\leqq \frac{1}{c}$. Our aim is to show that there exists a space $F\subseteqq \bigcap_{\nu \geqq \nu_0}  L_\nu^2(\R;H) $, which is dense in $L_\nu^2(\R;H)$ for all $\nu\geqq \nu_0$ and such that $\tilde S_\nu f=\tilde S_{\nu'}f$ for all $\nu,\nu'\geqq \nu_0$ and $f\in  F$. Having done so, we get the assertion from $\overline{\tilde S_\nu}=S_\nu$, $\nu\geqq \nu_0$, and Lemma \ref{le:indep_of_nu_for_bounded_mappings} (Use the Theorems \ref{thm:genSolth} and \ref{thm:causal_cont_causal} to get that $\tilde S_\nu|_F$ is bounded evolutionary at $\nu_0$ and (strongly) causal).


Let $\nu\geqq \nu_0$, $\s B_\nu\coloneqq \overline{\partial_0 \s M_\nu+\s N_\nu+\s A_\nu}$ and denote the domain of $\s B_\nu$ endowed with the graph norm by $D_{\s B_\nu}$. By Lemma \ref{le:nice_density}, we deduce that $\s F\coloneqq \lin \bigcup_{\delta>0}\left(1+\delta\partial_0\right)^{-1}[D(\s A)]$ is dense in $D_{\s B_\nu}$.  The mapping 
\[
   \iota \colon D_{\s B_\nu} \to L_\nu^2(\R;H), f\mapsto \s B_\nu f
\]
is continuous and onto, by Theorem \ref{thm:genSolth}. In particular, we have $\iota[\overline{G}]\subseteqq \overline{\iota[G]}$ for all $G\subseteqq D(\s B_\nu)$. Consequently, the set $\iota[\s F]$ is dense in $L_\nu^2(\R;H)$. Moreover, the set $\iota[\s F]$ lies in all $L_\nu^2(\R;H)$ for sufficiently large $\nu$. Moreover, $\s B_\nu f$ for all $f\in \s F$ is independent of $\nu$. Indeed,
\begin{align*}
   \s B_\nu f &= \overline{\left(\partial_0 \s M_\nu+\s N_\nu+\s A_\nu\right)}f \\
              &= \left(\partial_0 \s M_\nu+\s N_\nu+\s A_\nu\right)f \\
              &= \partial_0 \s M_\nu f+\s N_\nu f+\s A_\nu f \\
              &= \partial_0 \s M_\nu f+\s N_\nu f+\s A f,
\end{align*}
 where we have used that $(1+\delta\partial_0)^{-1}\s A_\nu \subseteqq \s A_\nu (1+\delta\partial_0)^{-1}$ and thus $\s F\subseteqq D(\s A)$ for all $\delta>0$. By Lemma \ref{le:indep_of_nu_for_bounded_mappings}, we get that $\s M_\nu f$ and $\s N_\nu f$ do not depend depend on $\nu$ and hence so does $\partial_{0}\s M_\nu f$. Thus, $F\coloneqq \iota[\s F]$ is the space we desired for. 
\end{proof}

\section{Proof of Theorem \ref{thm:genSolth}}\label{sec:proof_of_solth}

From now on, we assume the situation in Theorem \ref{thm:genSolth} as our \emph{standing hypothesis}. We need several preparations to give a proof for Theorem \ref{thm:genSolth}. Throughout, let $\s B\coloneqq \overline{\partial_0 \s M+\s N+\s A}$ and denote the domain of $\s B$ endowed with the graph norm by $D_{\s B}$. We denote the commutator of two operators $A,B$ by $[A,B]=AB-BA$ with its natural domain. Recall that in Hypothesis \ref{hyp:mat_law}, we assumed the existence of a continuous $M$ such that $[\s M,\partial_0]= M$ on $D(\partial_0)$. 
\begin{Le}\label{le:commus_of_a} Under the standing hypothesis, for every $\eps>0$, we have:
$$\overline{[(1+\eps\partial_0)^{-1},\partial_0 \s M]}u=\eps\partial_0 (1+\eps\partial_0)^{-1}M(1+\eps\partial_0)^{-1}u$$ for all $u\in L_\nu^2(\R;H)$.
\end{Le}
\begin{proof}
  Let $u\in D(\partial_0)$. We compute that
\begin{align*}
  [(1+\eps\partial_0)^{-1},\partial_0 \s M]u&= (1+\eps\partial_0)^{-1}\partial_0 \s M u - \partial_0 \s M (1+\eps\partial_0)^{-1}u\\
            & =(1+\eps\partial_0)^{-1}\left(\partial_0 \s M(1+\eps\partial_0)-(1+\eps\partial_0)\partial_0 \s M \right)(1+\eps\partial_0)^{-1}u \\
            & =\eps\partial_0(1+\eps\partial_0)^{-1} M(1+\eps\partial_0)^{-1}u.\qedhere
\end{align*}
\end{proof}
\begin{Le}\label{le:d(partial_0)cap d(A) is weakly dense}
  Under the standing hypothesis, we have for $\eps>0$ and $u\in D( \s B)$ that $(1+\eps\partial_0)^{-1}u\in D(\partial_0)\cap D(\s A)$. Moreover, the formula
\[
   \left(1+\eps\partial_0\right)^{-1}\s B u = \s B\left(1+\eps\partial_0\right)^{-1}u + \eps\partial_0(1+\eps\partial_0)^{-1}M\left(1+\eps\partial_0\right)^{-1}u 
  + \left[ \left(1+\eps\partial_0\right)^{-1},\s N\right] u 
\]
holds. In particular, we have 
\[
  \eps\partial_0(1+\eps\partial_0)^{-1}M\left(1+\eps\partial_0\right)^{-1}u \rightharpoonup 0
\]
and
\[
 \s B\left(1+\eps\partial_0\right)^{-1}u \rightharpoonup \s Bu
\]
as $\eps\to 0$ for all $u\in D(\s B)$ with weak convergence in $L_\nu^2(\R;H)$.
\end{Le}
\begin{proof}
 The computation can be made precise in $H_{-1}(\partial_0)\cap H_{-1}(\abs{\s A}+i)$\footnote{Here we use the concept of Sobolev lattices as introduced in \cite[Chapter 2]{Picard}. We briefly recall that for a densely defined closed linear operator $B\colon D(B)\subseteqq H\to H$ in some Hilbert space $H$ with $0\in \rho(B)$ we define $H_{-1}(B)$ to be the completion of $(H,\abs{B^{-1}\cdot}_H)$. It turns out that $H_{-1}(B)\cong D_{B^*}^*$ and that it is possible to continuously extend $B$ as a (unitary) mapping from $H$ to $H_{-1}(B)$. We note that if $B$ is normal, i.e., it commutes with its adjoint, then $H_{-1}(B)\cong D_B^*$.} Thus,
the commutator relation is a consequence of Lemma \ref{le:commus_of_a}. The convergence result relies on weak compactness in $L_\nu^2(\R;H)$, similar to the argument in \cite[Lemma 2.9]{RainerPicard2013}.
\end{proof}

\begin{Le}\label{le:ineq:forall}
  Under the standing hypothesis, we assume in addition that we are given a continuous $G \colon L_\nu^2(\R;H) \to \R$ and a bounded and measurable $\psi\colon\R\to\mathbb{C}$. Assume that for all $u\in D(\partial_0)\cap D(\s A)$ we have the inequality
\[
   \Re \langle \left(\partial_0 \s M +\s N +\s A \right) u,\psi(m_0)u\rangle \geqq G(u).
\]
Then the same inequality holds for $u\in D(\s B)$.
\end{Le}
\begin{proof}
  Let $u\in D(\s B)$. Using Lemma \ref{le:d(partial_0)cap d(A) is weakly dense}, we compute
\begin{align*}
   & \Re \langle \s B u,\psi(m_0)u\rangle \\
   & =\lim_{\eps\to 0}\Re \langle  \left(\partial_0 \s M +\s N +\s A \right)\left(1+\eps\partial_0\right)^{-1} u,\psi(m_0)\left(1+\eps\partial_0\right)^{-1}u\rangle \\
   & \geqq \lim_{\eps\to 0}G\left(\left(1+\eps\partial_0\right)^{-1} u\right) \\
   & = G(u) \qedhere
\end{align*}
\end{proof}
We come to the proof of Theorem \ref{thm:genSolth}:
\begin{proof}[Proof of Theorem \ref{thm:genSolth}]
  At first, we note that \[\Re \langle \s B u ,u \rangle \geqq c\langle u,u\rangle\] for all $u\in D(\s B)$ due to Lemma \ref{le:ineq:forall} yielding continuous invertibility once we have shown that $\s B$ is onto. For this, we let $\eps>0$, $f\in D(\s B^*)$ and show that $(1+\eps\partial_0^*)^{-1}[D(\s B^*)]\subseteqq D(\s B^*)$ and compute $\s B^*$. For $u\in D(\s B)$ we have
\begin{align*}
  &\langle \s B u,(1+\eps\partial_0^*)^{-1}f\rangle \\
  &= \langle \s B (1+\eps\partial_0)^{-1}u,f\rangle + \langle [\left(1+\eps\partial_0\right)^{-1},\partial_0 \s M]u,f\rangle + \langle [(1+\eps\partial_0)^{-1},\s N]u,f\rangle\\
  &=  \langle  u,(1+\eps\partial_0^*)^{-1}\s B^*f\rangle + \langle u,[\left(1+\eps\partial_0\right)^{-1},\partial_0 \s M]^*f\rangle + \langle u,[(1+\eps\partial_0)^{-1},\s N]^*f\rangle,
\end{align*}
proving that $(1+\eps\partial_0^*)^{-1}f \in D(\s B^*)$ and 
\begin{equation}
\label{adjoint}
  \s B^*\left(1+\eps\partial_0^*\right)^{-1}f = \left(1+\eps\partial_0^*\right)^{-1}\s B^*f +[\left(1+\eps\partial_0\right)^{-1},\partial_0 \s M]^*f+[(1+\eps\partial_0)^{-1},\s N]^*f. 
\end{equation}
Further, for $u\in D(\s A)\cap D(\partial_0)$ we compute 
\[
 \langle \left(\partial_0\s M + \s N + \s A\right) u,(1+\eps\partial_0^*)^{-1}f\rangle = \langle \s A u,(1+\eps\partial_0^*)^{-1}f\rangle +\langle u,\left(\partial_0\s M + \s N\right)^*(1+\eps\partial_0^*)^{-1}f\rangle.
\]
Using that $D(\s A)\cap D(\partial_0)$ is dense in $D(\s A)$ with respect to the graph norm of $\s A$, we get that $\left(1+\eps\partial_0^*\right)^{-1}f\in D(\s A^*)$ and 
\begin{equation}\label{adjoint2}
   \s B^*\left(1+\eps\partial_0^*\right)^{-1}f = \left(\left(\partial_0\s M + \s N\right)^*+\s A^*\right)(1+\eps\partial_0^*)^{-1}f. 
\end{equation}
Now, since $\left(1+\eps\partial_0^*\right)^{-1}\to 1$ strongly as $\eps\to 0$ and computing the adjoint is continuous with respect to the weak operator topology, we infer with the help of equation \eqref{adjoint} and Lemma \ref{le:d(partial_0)cap d(A) is weakly dense} that
\[
 \s B^*f=\textnormal{w-}\lim_{\eps\to0} \left(1+\eps\partial_0^*\right)^{-1}\s B^*f = \textnormal{w-}\lim_{\eps\to0} \s B^*\left(1+\eps\partial_0^*\right)^{-1}f,
\]
with limits taken in $L_\nu^2(\R;H)\subseteqq H_{-1}(\partial_0^*)\cap H_{-1}(\s A^*)$.
Moreover, 
\[
 \left(\left(\partial_0\s M + \s N\right)^*+\s A^*\right)(1+\eps\partial_0^*)^{-1}f\to \left(\left(\partial_0\s M + \s N\right)^*+\s A^*\right)f
\]
in $H_{-1}(\partial_0^*)\cap H_{-1}(\abs{\s A^*}+i)$. Thus, 
\[
 \s B^*\subseteqq \left(\partial_0\s M + \s N\right)^*+\s A^*,
\]
where the latter operator is considered with maximal domain in $L_\nu^2(\R;H)$, which equals $\{\phi\in L_\nu^2(\R;H); -\partial_0\s M^*\phi +\s A^*\phi\in L_\nu^2(\R;H)\}$. Observing that for all $u\in D(\s B^*)$:
\begin{align*}
  \Re \langle \s B^* u,u\rangle & = \lim_{\eps\to 0} \Re \langle \s B^*\left(1+\eps\partial_0^*\right)^{-1} u,\left(1+\eps\partial_0^*\right)^{-1}u\rangle \\
                                & = \lim_{\eps\to 0} \Re \langle \left(\left(\partial_0\s M + \s N\right)^*+\s A^*\right)\left(1+\eps\partial_0^*\right)^{-1} u,\left(1+\eps\partial_0^*\right)^{-1}u\rangle \\
   & \geqq \lim_{\eps\to 0} c\langle \left(1+\eps\partial_0^*\right)^{-1} u,\left(1+\eps\partial_0^*\right)^{-1}u\rangle =c\langle u,u\rangle,
\end{align*}
we deduce that $\s B^*$ is one-to-one, which, in turn, shows that $\s B$ is onto. 
\end{proof}

We conclude with a proof of Lemma \ref{le:nice_density}:
\begin{proof}[Proof of Lemma \ref{le:nice_density}]
 Let $u\in  D(\s B)$. For $\delta>0$, $\psi\in L_\nu^2(\R;H)$, $v\in D(\partial_0)\cap D(\s A)$ and $f\in F$, we compute
\begin{align*}
  & \abs{\langle \s B\left(u- \left(1+\delta\partial_0\right)^{-1}f\right),\psi\rangle} \\
   &\leqq \abs{\langle \s B \left(u- v\right),\psi\rangle}+ \abs{\langle \s B \left(v- \left(1+\delta\partial_0\right)^{-1}v\right),\psi\rangle}\\ &\quad +\abs{\langle \s B\left(\left(1+\delta\partial_0\right)^{-1}v- \left(1+\delta\partial_0\right)^{-1}f\right),\psi\rangle} \\
   &\leqq \abs{\langle \s B \left(u- v\right),\psi\rangle}+\abs{\langle \s B \left(v- \left(1+\delta\partial_0\right)^{-1}v\right),\psi\rangle}+\abs{\langle \partial_0\s M \left(1+\delta\partial_0\right)^{-1}\left(v- f\right),\psi\rangle} \\ &\quad +\abs{\langle \s N \left(1+\delta\partial_0\right)^{-1}\left(v- f\right),\psi\rangle}+\abs{\langle \s A  \left(1+\delta\partial_0\right)^{-1}\left(v- f\right),\psi\rangle}\\
   &\leqq \abs{\s B \left(u- v\right)}\abs{\psi}+\abs{\langle \s B \left(v- \left(1+\delta\partial_0\right)^{-1}v\right),\psi\rangle}\\ &\quad + \abs{\langle \s M \partial_0 \left(1+\delta\partial_0\right)^{-1}\left(v- f\right),\psi\rangle}+\abs{\langle M \left(1+\delta\partial_0\right)^{-1}\left(v- f\right),\psi\rangle} \\ &\quad +\abs{\langle \s N \left(1+\delta\partial_0\right)^{-1}\left(v- f\right),\psi\rangle} +\abs{\langle \left(1+\delta\partial_0\right)^{-1} \s A\left(v- f\right),\psi\rangle}\\ 
   &\leqq \abs{\s B \left(u- v\right)}\abs{\psi}+\abs{\langle \s B \left(v- \left(1+\delta\partial_0\right)^{-1}v\right),\psi\rangle}\\ &\quad + \frac{2}{\delta}\Abs{\s M}\abs{v- f}\abs{\psi}+\Abs{M}\abs{v- f}\abs{\psi} +\Abs{\s N}\abs{v- f}\abs{\psi}+\abs{\s A\left(v- f\right)}\abs{\psi}.
\end{align*}
By appropriately choosing $v$ such that $\abs{\s B (u-v)}$ is small (see Lemma \ref{le:d(partial_0)cap d(A) is weakly dense}) and afterwards choosing $\delta>0$ such that $\abs{\langle \s B \left(v- \left(1+\delta\partial_0\right)^{-1}v\right),\psi\rangle}$ is small, we find $f$ such that the remaining terms can be made small. Since $\s B$ is continuously invertible, for a suitable choice of sequences $(\delta_n)_n$ and $(f_n)_n$ we get that both $\left(u- \left(1+\delta_n\partial_0\right)^{-1}f_n\right)$ and $\s B\left(u- \left(1+\delta_n\partial_0\right)^{-1}f_n\right)$ are weakly convergent, which shows that $\bigcup_{\delta>0}\left(1+\delta\partial_0\right)^{-1}[F]$ is weakly dense in $D_{\s B}$ and hence $\s F$ is dense in $D_{\s B}$.
%
\end{proof}

 \bibliographystyle{plain}

\end{document}